\documentclass{rspublic}
\usepackage[latin1]{inputenc}
\usepackage{epsfig}
\usepackage{amssymb}
\usepackage{amsmath}

\begin{document}

\title{The first
digit frequencies of primes and Riemann zeta zeros tend to
uniformity following a size-dependent generalized Benford's law}

\author{Bartolo Luque, Lucas Lacasa}
\affiliation{Departamento de Matem{\'a}tica Aplicada y Estad{\'i}stica\\
ETSI Aeron{\'a}uticos, Universidad Polit{\'e}cnica de Madrid\\
28040 Madrid, Spain.}

\label{firstpage}

%% The \maketitle command is necessary to build the title page.
\maketitle

%%%%%%%%%%%%%%%%%%%%%%%%%%%%%%%%%%%%%%%%%%%%%%%%%%%%%%%%%%%%%%%%

\begin{abstract}{first significant digit, Benford's law, prime number,
pattern, zeta function.} Prime numbers seem to distribute among the
natural numbers with no other law than that of chance, however its
global distribution presents a quite remarkable smoothness. Such
interplay between randomness and regularity has motivated scientists
of all ages to search for local and global patterns in this
distribution that eventually could shed light into the ultimate
nature of primes. In this work we show that a generalization of the
well known first-digit Benford's law, which addresses the rate of
appearance of a given leading digit $d$ in data sets, describes with
astonishing precision the statistical distribution of leading digits
in the prime numbers sequence. Moreover, a reciprocal version of
this pattern also takes place in the sequence of the nontrivial
Riemann zeta zeros. We prove that the prime number theorem is, in
the last analysis, the responsible of these patterns. Some new
relations concerning the prime numbers distribution are also
deduced, including a new approximation to the counting function
$\pi(n)$. Furthermore, some relations concerning the statistical
conformance to this generalized Benford's law are derived. Some
applications are finally discussed.
\end{abstract}

\section{Introduction}
 \noindent The individual location of prime numbers within the integers
seems to be random, however its global distribution exhibits a
remarkable regularity (Zagier 1977). Certainly, this tenseness
between local randomness and global order has lead the distribution
of primes to be, since antiquity, a fascinating problem for
mathematicians (Dickson 2005) and more recently for physicists (see
for instance Berry \emph{et al.} 1999, Kriecherbauer \emph{et al}
2001, Watkins). The Prime Number Theorem, that addresses the global
smoothness of the counting function $\pi(n)$ providing the number of
primes less or equal to integer $n$, was the first hint of such
regularity (Tenenbaum 2000). Some other prime patterns have been
advanced so far, from the visual Ulam spiral (Stein \emph{et al}
1964) to the arithmetic progression of primes (Green \emph{et al}
2008), while some others remain conjectures, like the global gap
distribution between primes or the twin primes distribution
(Tenenbaum 2000), enhancing the mysterious interplay between
apparent randomness and hidden regularity. There are indeed many
open problems to be solved, and the prime number distribution is yet
to be understood (see for instance Guy 2004, Ribenboim 2004,
Caldwell). For instance, there exist deep connections between the
prime number sequence and the nontrivial zeros of the Riemann zeta
function (Watkins, Edwards 1964). The celebrated Riemann Hypothesis,
one of the most important open problem in mathematics, states that
the nontrivial zeros of the complex-valued Riemann zeta function
$\zeta(s)=\sum_{n=1}^\infty 1/n^s$ are all complex numbers with real
part $1/2$, the location of these being
intimately connected with the prime number distribution (Edwards 1964, Chernoff 2000).\\
\noindent Here we address the statistics of the first significant or
leading digit of both the sequences of primes and the sequence of
Riemann nontrivial zeta zeros. We will show that while the first
digit distribution is asymptotically uniform in both sequences (that
is to say, integers $1,...,9$ tend to be equally likely first digits
in both sequences when we take into account the infinite amount of
them), this asymptotic uniformity is reached in a very precise
trend, namely by following a size-dependent Generalized Benford's
law, what constitutes an as yet unnoticed pattern in both sequences.
The rest of the paper is organized as follows: in section 2 we
introduce the most celebrated first digit distribution: the
Benford's law. In section 3 we introduce a generalization of the
Benford's law, and we show that both the prime numbers and Riemann
zeta zeros sequences follow what we call a size-dependent
Generalized Benford's law, introducing two unnoticed patterns of
statistical regularity. In section 4 we point out that the mean
local density of both sequences is the responsible of these latter
patterns. We provide both statistical arguments (statistical
conformance between distributions) and analytical developments
(asymptotic expansions) that support our claim. In section 5 we
conclude and discuss on the possible applications.

\section{Benford's law}
The leading digit of a number stands for its non-zero leftmost
digit. For instance, the leading digits of the prime $7703$ and the
zeta zero $21.022...$ are $7$ and $2$ respectively. The most
celebrated leading digit distribution is the so called Benford's law
(Hill 1996), after physicist Frank Benford (1938) who empirically
found that in many disparate natural data sets and mathematical
sequences, the leading digit $d$ wasn't uniformly distributed as
might be expected, but instead had a biased probability of
appearance
\begin{equation}
P(d)=\log_{10}(1 + 1/d), \label{leydebenford}
\end{equation}
 where $d = 1,2,\dots, 9$. While this empirical law was indeed firstly discovered
by astronomer Simon Newcomb (1881), it is popularly known as the
Benford's law or alternatively as the Law of Anomalous Numbers.
Several disparate data sets such as stock prices, freezing points of
chemical compounds or physical constants exhibit this pattern at
least empirically. While originally being only a curious pattern
(Raimi 1976), practical implications began to emerge in the 1960s in
the design of efficient computers (see for instance Knuth 1967). In
recent years goodness-of-fit test against Benford's law has been
used to detect possible fraudulent financial data, by analyzing the
deviations of accounting data, corporation incomes, tax returns or
scientific experimental data to theoretical Benford predictions
(Nigrini 2000). Indeed, digit pattern analysis can produce valuable
findings not revealed by a mere glance, as is the
case of recent election results (Mebane 2006, Nigrini 2000).\\

\begin{figure}[h]
\centering
\includegraphics[width=0.8\textwidth]{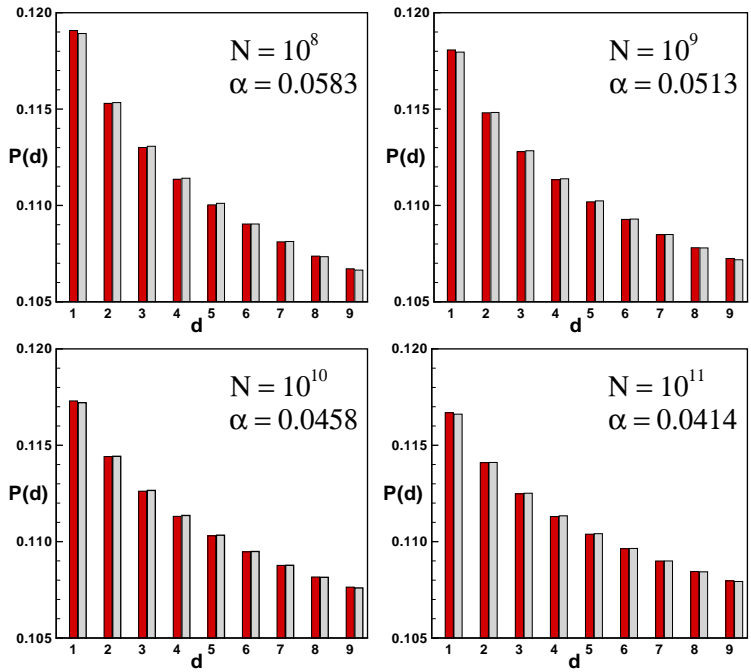}
\hspace{2cm} \caption{\small{Leading digit histogram of the prime
number sequence. Each plot represents, for the set of prime numbers
comprised in the interval $[1,N]$, the relative frequency of the
leading digit $d$ (red bars). Sample sizes are: $5761455$ primes for
$N=10^8$, $50847534$ primes for $N=10^9$, $455052511$ primes for
$N=10^{10}$ and $4118054813$ primes for $N=10^{11}$. Grey bars
represent the fitting to a Generalized Benford distribution (eq.
\ref{benfordgen}) with a given exponent $\alpha(N)$. }}
\label{4ajustes}
\end{figure}

\noindent Many mathematical sequences such as $(n^n)_{n \in
\mathbb{N}}$ and $(n!)_{n \in \mathbb{N}}$ (Benford 1938), binomial
arrays $(^n_k)$ (Diaconis  1977), geometric sequences or sequences
generated by recurrence relations (Raimi 1976, Miller \emph{et al.}
2006) to cite a few are proved to be Benford. One may thus wonder if
this is the case for the primes. In figure \ref{4ajustes} we have
plotted the leading digit $d$ rate of appearance for the prime
numbers placed in the interval $[1,N]$ (red bars), for different
sizes $N$. Note that intervals $[1,N]$ have been chosen such that
$N=10^D,\ D\in\mathbb{N}$ in order to assure an unbiased sample
where all possible first digits are equiprobable a priori (see the
appendix for further details). Benford's law states that the first
digit of a series data extracted at random is $1$ with a frequency
of $30.1\%$, and is $9$ only about $4.6\%$. Note in figure
\ref{4ajustes} that primes seem however to approximate uniformity in
its first digit. Indeed, the more we increase the interval under
study, the more we approach uniformity (in the sense that all
integers $1,...,9$ tend to be equally likely as a first digit). As a
matter of fact, Diaconis (1977) proved that primes are not Benford
distributed as long as their first significant digit is
asymptotically uniformly distributed. A question arises
straightforwardly: how does the prime sequence reach this uniform
behavior in the infinite limit? Is there any pattern on its trend
towards uniformity, or on the contrary, does the first digit
distribution lacks any structure for finite sets?

\section{Generalized Benford's law: the pattern}
Several mathematical insights of the Benford's law have been also
advanced so far (Hill 1995a, Pinkham 1961, Raimi 1976, Miller
\emph{et al.} 2006), and Hill (1995b) proved a Central Limit-like
Theorem which states that random entries picked from random
distributions form a sequence whose first digit distribution
converges to the Benford's law, explaining thereby its ubiquity.
This law has been for a long time practically the only distribution
that could explain the presence of skewed first digit frequencies in
generic data sets. Recently Pietronero \emph{et al.} (2001) proposed
a generalization of Benford's law based in multiplicative processes
(see also Nigrini \emph{et al.} 2007). It is well known that a
stochastic process with probability density $1/x$ generates data
which are Benford, therefore series generated by power law
distributions $P(x) \sim x^{-\alpha}$ with $\alpha \neq1$, would
have a first digit distribution that follow a so-called Generalized
Benford's law (GBL):
\begin{equation}
P(d)= C\int_d^{d+1}x^{-\alpha}dx=\frac{1}{10^{1-\alpha}-1} \bigg
[(d+1)^{1-\alpha}-d^{1-\alpha}\bigg],\label{benfordgen}
\end{equation}
where the prefactor is fixed for normalization to hold and $\alpha$
is the exponent of the original power law distribution (for
$\alpha=1$, the GBL reduces to the Benford's law).
\subsection{\textbf{The pattern in primes}}
\begin{figure}[h]
\centering
\includegraphics[width=0.8\textwidth]{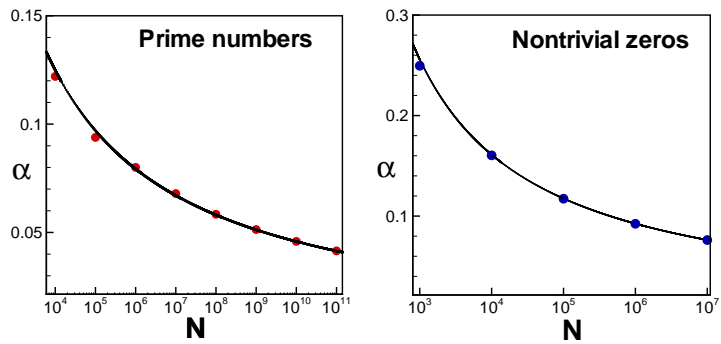}
\hspace{2cm} \caption{\small{Size dependent parameter $\alpha(N)$.
Left: Red dots represent the exponent $\alpha(N)$ for which the
first significant digit of prime number sequence fits a Generalized
Benford Law in the interval $[1,N]$. The black line corresponds to
the fitting, using a least squares method, $\alpha(N)=1/(\log
N-1.10)$. Right: same analysis as for the left figure, but for the
Riemann nontrivial zeta zeros sequence. The best fitting is
$\alpha(N)=1/(\log N-2.92)$.}} \label{alfas}
\end{figure}

\noindent Although Diaconis showed that the leading digit of primes
distributes uniformly in the infinite limit, there exist a clear
bias from uniformity for finite sets (see figure \ref{4ajustes}). In
this figure we have also plotted (grey bars) the fitting to a GBL.
Note that in each of the four intervals, there is a particular value
of exponent $\alpha$ for which an excellent agreement holds (see the
appendix for fitting methods and statistical tests). More
specifically, given an interval $[1,N]$, there exists a particular
value $\alpha(N)$ for which a GBL fits with extremely good accuracy
the first digit distribution of the primes appearing in that
interval. Interestingly, the value of the fitting parameter $\alpha$
decreases as the interval, hence the number of primes, increases in
a very particular way. In the left part of figure \ref{alfas} we
have plotted this size dependence, showing that a functional
relation between $\alpha$ and $N$ takes place:
\begin{equation}
\alpha(N)= \frac{1}{\log N - a}, \label{alfa_n}
\end{equation}
where $a=1.10\pm0.05$ for large values of $N$. Notice that
$\lim_{N\rightarrow\infty}\alpha(N)=0$, and in this situation this
size-dependent GBL reduces to the uniform distribution, in
consistency with previous theory (Diaconis 1977). Despite the local
randomness of the prime numbers sequence, it seems that its first
digit distribution converges smoothly to uniformity in a very
precise trend: as a GBL with a size dependent exponent
$\alpha(N)$.\\
\subsubsection{GBL Extension} 
\begin{figure}[h]
\centering
\includegraphics[width=0.8\textwidth]{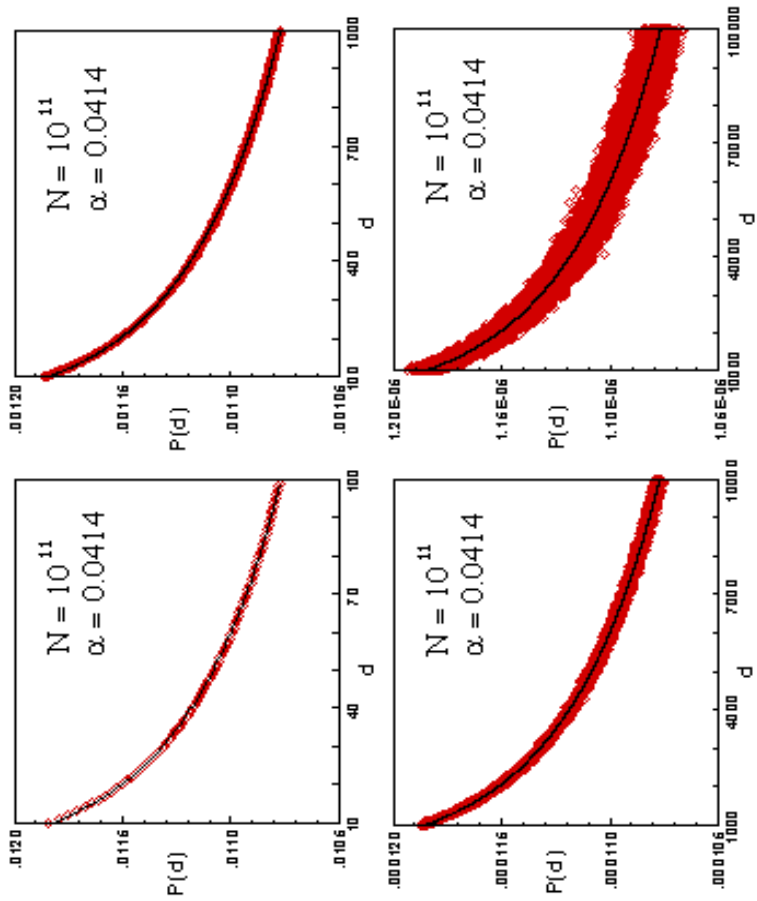}
\hspace{2cm} \caption{\small{Extension of GBL to the $k$ first
significant digits. In this figure we represent the fitting of an
extended GBL following eq. \ref{kdigits} (black line) to the first
two significant digits relative frequencies (up-left), first three
significant digits relative frequencies (up-right), first four
significant digits relative frequencies (down-left) and first five
significant digits relative frequencies (down-right) of the
$4118054813$ primes appearing in the interval $[1,10^{11}]$ (red
dots).}}\label{extension}
\end{figure}

At this point an extension of the GBL
to include, not only the first significative digit, but the first
$k$ significative ones can be done (Hill 1995b). Given a number $n$,
we can consider its $k$ first significative digits $d_1,d_2 \dots,
d_k$ through its decimal representation: $D=\sum_{i=1}^k{d_i
10^{k-i}}$, where $d_1 \in \{1,\dots, 9 \}$ and $d_i \in
\{0,1,\dots, 9 \}$ for $i \geq 2$. Hence, the \emph{extended} GBL
providing the probability of starting with number $D$ is
\begin{eqnarray}
&&P(d_1, d_2, \dots, d_k)=P(D)\nonumber \\
&&= \frac{1}{(10^k)^{1-\alpha}-10^{k-1}}\bigg
[(D+1)^{1-\alpha}-D^{1-\alpha} \bigg]. \label{kdigits}
\end{eqnarray}
Figure \ref{extension} represents the fitting of the $4118054813$
primes appearing in the interval $[1,10^{11}]$ to an \emph{extended}
GBL for $k=2,3,4$ and $5$: interestingly, the pattern still holds.\\

\subsection{\textbf{The `mirror' pattern in the Riemann zeta zeros sequence}}
\noindent Once the pattern has been put forward in the case of the
prime number sequence, we may wonder if a similar behavior holds for
the sequence of nontrivial Riemann zeta zeros (zeros sequence from
now on). This sequence is composed by the imaginary part of the
nontrivial zeros (actually only those with positive imaginary part
are taken into account by symmetry reasons) of $\zeta(s)$. While
this sequence is not Benford distributed in the light of a theorem
by Rademacher-Hlawka (1984) that proves that it is
\emph{asymptotically} uniform, will it follow a size-dependent GBL
as in the case of the primes?

\begin{figure}[h]
\centering
\includegraphics[width=0.8\textwidth]{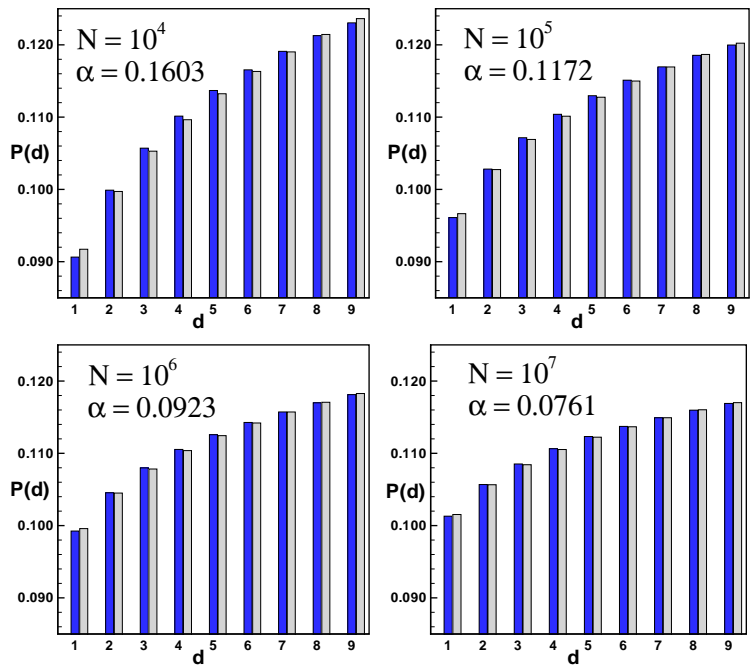}
\hspace{2cm} \caption{\small{Leading digit histogram of the
nontrivial Riemann zeta zeros sequence. Each plot represents, for
the sequence of Riemann zeta zeros comprised in the interval
$[1,N]$, the observed relative frequency of leading digit $d$ (blue
bars). Sample sizes are: $10142$ zeros for $N=10^4$, $138069$ zeros
for $N=10^5$, $1747146$ zeros for $N=10^6$ and $21136126$ zeros for
$N=10^7$. Grey bars represent the fitting to a GBL following
equation \ref{gbl2} with a given exponent $\alpha(N)$. }}
\label{4ajusteszero}
\end{figure}

\begin{figure}[h]
\centering
\includegraphics[width=0.8\textwidth]{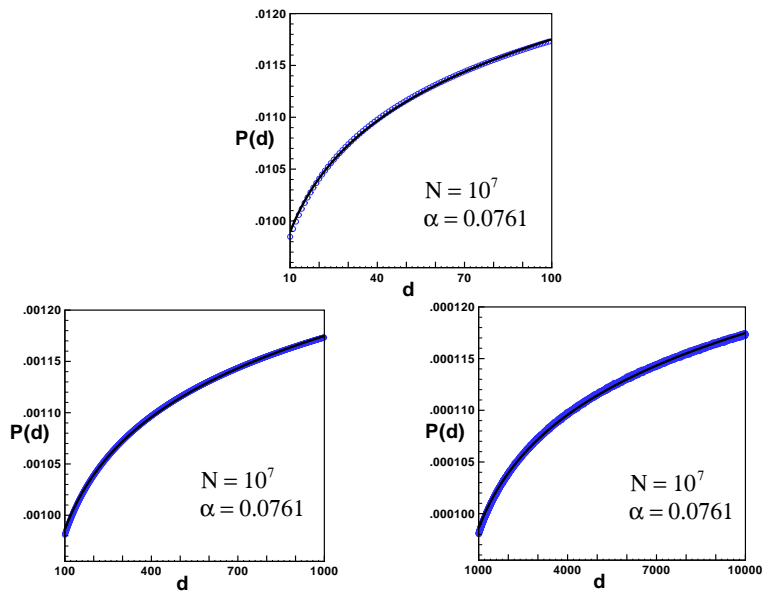}
\hspace{2cm} \caption{\small{Extension of GBL to the $k$ first
significant digits. In this figure we represent the fitting of an
extended GBL following eq. \ref{kdigits2} (black line) to the first
two significant digits relative frequencies (up), first three
significant digits relative frequencies (down-left), and first four
significant digits relative frequencies (down-right) of the
$21136126$ zeros appearing in the interval $[1,10^7]$ (blue dots).
}}\label{extension_zeros}
\end{figure}

\noindent In figure \ref{4ajusteszero} we have plotted, in the
interval $[1,N]$ and for different values of $N$, the relative
frequencies of leading digit $d$ in the zeros sequence (blue bars),
and in grey bars a fitting to a GBL with density $x^{\alpha}$, i.e.:
\begin{equation}
P(d)= C\int_d^{d+1}x^{\alpha}dx=\frac{1}{10^{1+\alpha}-1}
\bigg[(d+1)^{1+\alpha}-d^{1+\alpha}\bigg] \label{gbl2}
\end{equation}
(this reciprocity is clarified later in the text). Note that a very
good agreement holds again for particular size-dependent values of
$\alpha$, and the same functional relation as equation \ref{alfa_n}
holds with $a=2.92\pm0.05$. As in the case of the primes, this size
dependent GBL tends to uniformity for $N\rightarrow\infty$, as it
should (Hlawka 1984). Moreover, the \emph{extended} version of
equation \ref{gbl2} for the $k$ first significative digits is
\begin{eqnarray}
&&P(d_1, d_2, \dots, d_k)=P(D)\nonumber \\
&&= \frac{1}{(10^k)^{1+\alpha}-10^{k-1}}\bigg
[(D+1)^{1+\alpha}-D^{1+\alpha} \bigg]. \label{kdigits2}
\end{eqnarray}
As can be seen in figure \ref{extension_zeros}, the pattern also
holds in this case.

\section{Explanation of the primes pattern}

\noindent Why do these two sequences exhibit this unexpected pattern
in the leading digit distribution? What is the responsible for it to
take place? While the prime number distribution is deterministic in
the sense that precise rules determine whether an integer is prime
or not, its apparent local randomness has suggested several
stochastic interpretations. In particular, Cram\'{e}r (1935, see
also Tenembaum 2000) defined the following model: assume that we
have a sequence of urns $U(n)$ where $n = 1,2,...$ and put black and
white balls in each urn such that the probability of drawing a white
ball in the $k^{th}$-urn goes like $1/\log k$. Then, in order to
generate a sequence of pseudo-random prime numbers we only need to
draw a ball from each urn: if the drawing from the $k^{th}$-urn is
white, then $k$ will be labeled as a pseudo-random prime. The prime
number sequence can indeed be understood as a concrete realization
of this stochastic process,
where the chance of a given integer $x$ to be prime is $1/\log x$.\\
We have repeated all statistical tests within the stochastic
Cram\'{e}r model, and have found that a statistical sample of
pseudo-random prime numbers in $[1,10^{11}]$ is also GBL distributed
and reproduce all statistical analysis previously found in the
actual primes (see the appendix for an in-depth analysis). This
result strongly suggests that a density $1/\log x$, which is nothing
but the mean local primes density by virtue of the prime number
theorem, is likely to be the responsible for the GBL pattern. In
what follows we will provide further statistical and analytical
arguments that support this fact.

%%%%%%%%%%%%%%%%%%%%%%SUPLEMENTAERIO/////////////
\subsection{\textbf{Statistical conformance of prime number distribution to
GBL}} \noindent Recently, it has been shown that disparate
distributions such as the Lognormal, the Weibull or the Exponential
distribution can generate standard Benford behavior (Leemis \emph{et
al.} 2000) for particular values of their parameters. In this sense,
a similar phenomenon could be taking place with GBL: can different
distributions generate GBL behavior? One should thus switch the
emphasis from the examination of data sets that obey GBL to
\emph{probability distributions} that do so, other than power laws.

\subsubsection{\textbf{$\chi^2$-test for conformance between distributions}}

\noindent The prime counting function $\pi(N)$ provides the number
of primes in the interval $[1,N]$ (Tenenbaum \emph{et al.} 2000) and
up to normalization, stands as the cumulative distribution function
of primes. While $\pi(N)$ is a stepped function, a nice asymptotic
approximation is the offset logarithmic integral:
\begin{equation}
\pi(N)\sim\int_2^N\frac{1}{\log x}dx= \textrm{Li}(N), \label{intlog}
\end{equation}
(one of the formulations of the Riemann hypothesis actually states
that $|\textrm{Li}(n)-\pi(n)|<c\sqrt n\log n$, for some constant $c$
(Edwards 1974)). We can interpret $1/\log x$ as an average prime
density and the lower bound of the integral is set to be $2$ for
singularity reasons. Following Leemis \emph{et al.} (2000), we can
calculate a chi-square goodness-of-fit test of the conformance
between the first digit distribution generated by $\textrm{Li}(N)$
and a GBL with exponent $\alpha(N)$. The test statistic is in this
case:
\begin{equation}
c=\sum_{d=1}^9\frac{[\textrm{Pr}(Y=d)-\textrm{Pr}(X=d)]^2}{\textrm{Pr}(X=d)},
\end{equation}
where $\textrm{Pr}(X)$ is the first digit probability (eq.
\ref{benfordgen}) for a GBL associated to a probability distribution
with exponent $\alpha(N)$ and $\textrm{Pr}(Y)$ is the tested
probability. In table \ref{chi2_pi_li} we have computed, fixed the
interval $[1,N]$, the chi-square statistic $c$ for two different
scenarios, namely the normalized logarithmic integral
$\textrm{Li}(n)/\textrm{Li}(N)$ and the normalized prime counting
function $\pi(n)/\pi(N)$, with $n \in [1,N]$. In both cases there is
a remarkable good agreement and we cannot reject the hypothesis that
primes are size-dependent GBL.
\begin{table}
\center
\begin{tabular}{|c|c|c|}
\hline $N$ & $c$ for &$c$ for\\

& $\pi(x)/\pi(N)$ & $\textrm{Li}(x)/\textrm{Li}(N)$\\
\hline
$10^3$ & $0.59\cdot10^{-2}$ & $0.58\cdot10^{-2}$ \\
$10^4$ & $0.86\cdot10^{-3}$ & $0.57\cdot10^{-3}$ \\
$10^5$ & $0.12\cdot10^{-3}$ & $0.13\cdot10^{-3}$\\
$10^6$ & $0.57\cdot10^{-4}$ & $0.61\cdot10^{-4}$ \\
$10^7$ & $0.32\cdot10^{-4}$ & $0.33\cdot10^{-4}$\\
$10^8$ & $0.17\cdot10^{-4}$ & $0.17\cdot10^{-4}$ \\
\hline
\end{tabular}
\caption{\small{Chi-square goodness-of-fit test $c$ of the
conformance between primes cumulative distributions ($\pi(x)/\pi(N)$
and $\textrm{Li}(x)/\textrm{Li}(N)$) and a GBL with exponent
$\alpha(N)$ (eq. \ref{alfa_n}) in the interval $[1,N]$. The null
hypothesis, prime number distribution obeys GBL, cannot be rejected
.}} \label{chi2_pi_li}
\end{table}

\subsubsection{\textbf{Conditions for conformance to GBL}}
\noindent Hill (1995b) wondered about which common distributions (or
mixtures thereof) satisfy Benford's law.  Leemis \emph{et al.}
(2000) tackled this problem and quantified the agreement to
Benford's law of several standard distributions. They concluded that
the ubiquity of Benford behavior could be related to the fact that
many distributions follow Benford's law for particular values of
their parameters. Here, following the philosophy of that work
(Leemis \emph{et al.} 2000), we will develop a mathematical
framework that provide conditions for
conformance to a GBL.\\

\noindent The probability density function of a discrete GB random
variable $Y$ is:
\begin{equation}
f_Y(y)=\textrm{Pr}(Y=y)=\frac{1}{10^{1-\alpha}-1}[(y+1)^{1-\alpha}-y^{1-\alpha}],\
y=1,2,...,9. \label{BGdemo}
\end{equation}
The associated cumulative distribution function is therefore:
\begin{equation}
F_Y(y)=\textrm{Pr}(Y\leq
y)=\frac{1}{10^{1-\alpha}-1}[(y+1)^{1-\alpha}-1],\  y=1,2,...,9.
\label{BGacumulada}
\end{equation}
How can we prove that a random variable $T$ extracted from a
probability density $f_T(t)=\textrm{Pr}(t)$ has an associated
(discrete) random variable $Y$ that follows equation \ref{BGdemo}?
We can readily find a relation between both random variables.
Suppose without loss of generality that the random variable $T$ is
defined in the interval $[1,10^{D+1})$. Let the discrete random
variable $D$ fulfill:
\begin{equation}
10^D \leq T < 10^{D+1}
\end{equation}
This definition allows us to express the first significative digit
$Y$ in terms of $D$ and $T$:
\begin{equation}
Y=\lfloor T\cdot 10^{-D} \rfloor,
\end{equation}
where from now on the floor brackets stand for the integer part
function. Now, let $U$ be a random variable uniformly distributed in
$(0,1)$, $U \sim U(0,1)$. Then, inverting the cumulative
distribution function \ref{BGacumulada} we come to:
\begin{equation}
Y=\lfloor[(10^{1-\alpha}-1)\cdot U+1]^{\frac{1}{1-\alpha}}\rfloor.
\label{gene_BG}
\end{equation}

\noindent This latter relation is useful to generate a discrete GB
random variable $Y$ from a uniformly distributed one $U(0,1)$. Note
also that for $\alpha=0$, we have $Y=\lfloor9\cdot U+1\rfloor$, that
is, a first digit distribution which is uniform
$\textrm{Pr}(Y=y)=1/9,\ y=1,2,...,9$, as expected. Hence, every
discrete random variable $Y$ that distributes as a GB should fulfill
equation \ref{gene_BG}, and consequently if a random variable $T$
has an associated random variable $Y$, the following identity should
hold:
\begin{equation}
\lfloor T\cdot 10^{-D}\rfloor=\lfloor[(10^{1-\alpha}-1)\cdot
U+1]^{\frac{1}{1-\alpha}}\rfloor,
\end{equation}
and then,
\begin{equation}
Z= \frac{(T10^{-D})^{1-\alpha}-1}{10^{1-\alpha}-1} \sim U(0,1).
\end{equation}
In other words, in order the random variable $T$ to generate a GB,
the random variable $Z$ defined in the preceding transformation
should distribute as $U(0,1)$. The cumulative distribution function
of $Z$ is thus given by:
\begin{equation}
F_Z(z)= \sum_{d=0}^n\bigg\{\textrm{Pr}(10^d\leq T < 10^{d+1})\cdot
\textrm{Pr}\bigg(\frac{(T10^{-D})^{1-\alpha}-1}{10^{1-\alpha}-1}
\leq z | 10^d\leq T < 10^{d+1}\bigg)\bigg\}=z, \label{testBG1}
\end{equation}
that in terms of the cumulative distribution function of $T$ becomes
\begin{equation}
\sum_{d=0}^n\{{F_T(v10^d)-F_T(10^d)}\}=z, \label{testBG2}
\end{equation}
where $v\equiv[(10^{1-\alpha}-1)z+1]^{\frac{1}{1-\alpha}}$.\\

\noindent We may take now the power law density $x^{-\alpha}$
proposed by Pietronero \emph{et al.} (2001) in order to show that
this distribution exactly generates Generalized Benford behavior:
\begin{equation}
f_T(t)=\textrm{Pr}(t)=\frac{1-\alpha}{10^{(D+1)(1-\alpha)}-1}t^{-\alpha},
\ t \in [1,10^{D+1})
\end{equation}
Its cumulative distribution function will be:
\begin{equation}
F_T(t)=\frac{t^{1-\alpha}-1}{10^{(D+1)(1-\alpha)}-1},
\end{equation}
and thereby equation \ref{testBG2} reduces to:
\begin{equation}
\sum_{d=0}^D\{{F_T(v10^d)-F_T(10^d)}\}=\frac{z(10^{1-\alpha}-1)}{10^{(D+1)(1-\alpha)}-1}\sum_{d=0}^D{(10^{1-\alpha})^d}=z,
\end{equation}
as expected.

\subsubsection{\textbf{GBL holds for prime number distribution}} \noindent
While the preceding development is in itself interesting in order to
check for the conformance of several distributions to GBL, we will
restrict our analysis to the prime number cumulative distribution
function conveniently normalized in the interval $[1,10^D]$:
\begin{equation}
F_T(t)=\frac{\pi(t)}{\pi(10^{D+1})},\ t \in [1,10^{D+1})
\end{equation}
Note that previous analysis showed that
\begin{equation}
\alpha(10^{D+1})=\frac{1}{\ln{(10^{D+1})}-a},
\end{equation}
where $a\simeq1.1$. Since $\pi(t)$ is a stepped function that does
not possess a closed form, the relation \ref{testBG2} cannot be
analytically checked. However a numerical exploration can indicate
into which extent primes are conformal with GBL. Note that relation
\ref{testBG2} reduces in this case to
\begin{equation}
\sum_{d=0}^D\bigg\{\pi(v\cdot10^d)-\pi(10^d)\bigg\}=\pi(10^{D+1})z
\label{zprimos}
\end{equation}
where
$v\equiv[(10^{1-\alpha(10^{D+1})}-1)z+1]^{\frac{1}{1-\alpha(10^{D+1})}}$
and $z\in[0,1]$. Firstly, this latter relation is trivially
fulfilled for the extremal values $z=0$ and $z=1$. For other values
$z\in(0,1)$, we have numerically tested this equation for different
values of $D$, and have found that it is satisfied with negligible
error (we have performed a scatterplot of equation \ref{zprimos} and
have found a correlation coefficient
$r=1.0$).\\
The same numerical analysis has been performed for logarithmic Li.
integral. In this case the relation
\begin{equation}
\sum_{d=0}^D\bigg\{\textrm{Li}(v\cdot10^d)-\textrm{Li}(10^d)\bigg\}=\textrm{Li}(10^{D+1})z,
\label{zprimos2}
\end{equation}
is satisfied with similar remarkable results provided that we fix
$\textrm{Li}(1)\equiv0$ for singularity reasons.

%%%%%%%%%%%%%%%%%%%////////////

%We have confirmed (supporting information) that there exists
%statistical conformance between the prime number cumulative
%distribution $\pi(n)$ (conveniently normalized in $[1,N]$) and a GBL
%with exponent $\alpha(N)$. The same holds for the Eulerean
%logarithmic integral $\textrm{Li}(n)$, which constitutes an
%asymptotic approximation to $\pi(n)$:
%\begin{equation}
%\pi(N)\sim\textrm{Li}(N)=\int_2^N\frac{1}{\log x}dx \label{li1}
%\end{equation}
%(one of the formulations of the Riemann hypothesis actually states
%that $|\textrm{Li}(n)-\pi(n)|<c\sqrt n\log n$, for some constant $c$
%\cite{Edwards}).

%\noindent We can indeed prove (supporting information) that any
%cumulative distribution function $F_T(z)$ that evidences GBL
%behavior must fulfill
%\vspace*{-3.5mm}
%\begin{equation}
%\sum_{d=0}^n\{F_T(v10^d)-F_T(10^d)\}=z,
%\label{f}
%\end{equation}
%\vspace*{-0.5mm}
%\noindent where
%$v\equiv[(10^{1-\alpha}-1)z+1]^{\frac{1}{1-\alpha}}$ and
%$z\in[0,1]$. This latter relation is trivially fulfilled by power
%law densities $x^{-\alpha}$.

%In particular, eq. \ref{f} reduces for the prime cumulative
%distribution $\pi(n)$ to
%\begin{equation}
%\sum_{d=0}^D\left\{\pi(v\cdot10^d)-\pi(10^d)\right\}=\pi(10^{D+1})z,
%\label{zprimos}
%\end{equation}
%and
%\begin{equation}
%\sum_{d=0}^D\bigg\{\textrm{Li}(v\cdot10^d)-\textrm{Li}(10^d)\bigg\}=\textrm{Li}(10^{D+1})z,
%\label{zprimos2}
%\end{equation}
%for $\textrm{Li}(n)$. Equations \ref{zprimos} and \ref{zprimos2} are numerically
%fulfilled (supporting information).\\

\subsection{\textbf{Asymptotic expansions}}
\noindent Hitherto, we have provided statistical arguments that
indicate that other distributions than $x^{-\alpha}$ such as $1/\log
x$ can generate GBL behavior. In what follows we provide analytical
arguments that support this fact. \\$\textrm{Li}(N)$ possesses the
following asymptotic expansion
\begin{eqnarray}
\textrm{Li}(N)=\frac{N}{\log N}\bigg\{1+\frac{1}{\log
N}+\frac{2}{\log^2 N }+O\bigg(\frac{1}{\log^3
N}\bigg)\bigg\}.\label{li}
\end{eqnarray}
\noindent Now, a sequence whose first significant digit follows a
GBL has indeed a density that goes as $x^{-\alpha}$. One can
consequently derive from this latter density a function $L(N)$ that
provides the number of primes appearing in the interval $[1,N]$ as
it follows:
\begin{equation}
L(N)= e\alpha(N)\int_2^Nx^{-\alpha(N)}dx \label{L(N)}
\end{equation}
where the prefactor is fixed for $L(N)$ to fulfill the prime number
theorem and consequently
\begin{equation}
\lim_{N\rightarrow\infty}\frac{L(N)}{N/\log N}=1
\end{equation}

\begin{table}
\center
\begin{tabular}{|c|c|c|c|c|c|c|}
\hline
$N$ & $\pi(N)$ & $\textrm{Li}(N)$ & $N/\log N$ & $L(N)$ & $L(N)/\pi(N)$ \\
\hline
$10^2$ & $25$ & $30$ &    $22$ & $29$ & $0.85533$\\
$10^3$ & $168$ & $178$ &    $145$ & $172$ & $0.97595$\\
$10^4$ & $1229$ & $1246$ &    $1086$ & $1228$ & $1.00081$\\
$10^5$ & $9592$ & $9630$ &    $8686$ & $9558$ & $1.00352$\\
$10^6$ & $78492$ & $78628$ &  $72382$ & $78280$ & $1.00278$\\
$10^7$ & $664579$ & $664918$ &  $620421$ & $662958$ & $1.00244$\\
$10^8$ & $5761455$ & $5762209$ &    $5428681$ & $5749998$ & $1.00199$\\
$10^{9}$ & $50847534$ & $50849235$ & $48254942$ & $50767815$ & $1.00157$\\
$10^{10}$ & $455052511$ & $455055615$ & $434294482$ & $454484882$ & $1.00125$\\
$10^{20}$ & $2220819602560918840$ & $$ & $$ & $$ & $1.00027$\\
\hline

\end{tabular}
\caption{\small{Up to integer $N$, values of the prime counting
function $\pi(N)$, the approximation given by the logarithmic
integral $\textrm{Li}(N)$, $N/\log N$, the counting function $L(N)$
defined in eq. \ref{L(N)} and the ratio $L(N)/\pi(N)$.}}
\label{Table_counting}
\end{table}

(see table \ref{Table_counting} for a numerical exploration of this
new approximation to $\pi(N)$). Now, we can asymptotically expand
$L(N)$ as it follows
\begin{eqnarray}
&&L(N)=\frac{\alpha(N)e}{1-\alpha(N)}N^{1-\alpha(N)}\nonumber\\
&&=\frac{N}{\log N-(a+1)}\cdot \exp{\bigg(\frac{-a}{\log N-a}\bigg)}\nonumber\\
&&=\frac{N}{\log N}\bigg\{1+\frac{a+1}{\log N} +
\frac{(a+1)^2}{\log^2 N} +
O \bigg(\frac{1}{\log^3 N}\bigg)\bigg\}\cdot\nonumber\\
&&\cdot\bigg\{1-\frac{a}{\log N-a} + \frac{a^2}{(\log N-a)^2} +
O \bigg(\frac{1}{(\log N-a)^3}\bigg)\bigg\}\nonumber\\
&&=\frac{N}{\log N}\bigg\{1+\frac{1}{\log N}+\frac{1+a-a^2/2}{\log^2
N}+O \bigg(\frac{1}{\log^3 N}\bigg)\bigg\}.\label{lu}
\end{eqnarray}
Comparing equations \ref{li} and \ref{lu}, we conclude that
$\textrm{Li}(N)$ and $L(N)$ are compatible cumulative distributions
within an error
\begin{equation}
E(N)=\frac{N}{\log N}\bigg\{\frac{2}{\log^2 N
}-\frac{1+a-a^2/2}{\log^2 N}+O\bigg(\frac{1}{\log^3 N}\bigg)\bigg\}
\label{error}
\end{equation}
that is indeed minimum for $a=1$, in consistency with our previous
numerical results obtained for the fitting value of $a$ (eq.
\ref{alfa_n}). Hence, within that error we can conclude that primes
obey
a GBL with $\alpha(N)$ following equation \ref{alfa_n}: primes follow a size-dependent generalized Benford's law.\\

\section{Explanation of the pattern in the case of the Riemann zeta zeros sequence}
\noindent What about the Riemann zeros? Von Mangoldt proved (Edwards
1974) that on average, the number of nontrivial zeros
 $R(N)$ up to $N$ (zeros counting function) is
\begin{equation}
R(N)=\frac{N}{2\pi}\log{\bigg(\frac{N}{2\pi}\bigg)}-\frac{N}{2\pi} +
O(\log N).
\end{equation}
$R(N)$ is nothing but the cumulative distribution of the zeros (up
to normalization), which satisfies
\begin{equation}
R(N)\approx \frac{1}{2\pi} \int_{2}^N \log
\bigg(\frac{x}{2\pi}\bigg)dx. \label{zerosss}
\end{equation}
The nontrivial Riemann zeros average density is thus $\log(x/2\pi)$,
which is nothing but the reciprocal of the prime numbers mean local
density (see eq. \ref{intlog}). One can thus straightforwardly
deduce a power law approximation to the cumulative distribution of
the non trivial zeros similar to equation \ref{L(N)}:
\begin{equation}
R(N)\sim \frac{1}{2\pi e\alpha(N/2\pi)}\int_2^N
\bigg(\frac{x}{2\pi}\bigg)^{\alpha(N/2\pi)}dx.
\label{counting_zeros}
\end{equation}
We conclude that zeros are also GBL for $\alpha(N)$ satisfying the
following change of scale
\begin{equation}
\alpha(N/2\pi)=\frac{1}{\log(N/2\pi)-a}= \frac{1}{\log N -
(\log(2\pi)+ a)}.
\end{equation}
Hence, since $a\simeq1.1$ (equation \ref{error}) one should expect
for the constant $a$ associated to the zeros sequence the following
value: $\log(2\pi)+1.1 \approx 2.93$, in good agreement with our
previous numerical analysis.

\section{Discussion}
\noindent To conclude, we have unveiled a statistical pattern in the
prime numbers and the nontrivial Riemann zeta zeros sequences that
has surprisingly gone unnoticed until now. According to several
statistical and analytical arguments, we can conclude that the shape
of the mean local density of both sequences are the responsible of
these patterns. Along with this finding, some relations concerning
the statistical conformance of any given distribution to the
generalized Benford's law have also been derived.\\
Several applications and future work can be depicted: first, since
the Riemann zeros seem to have the same statistical properties as
the eigenvalues of a concrete type of random matrices called the
Gaussian Unitary Ensemble (Berry 1999, Bogomolny 2007), the relation
between GBL and random matrix theory should be investigated in-depth
(Miller \emph{et al.} 2005). Second, this finding may also apply to
several other sequences that, while not being strictly Benford
distributed, can be GBL, and in this sense much work recently
developed for Benford distributions (H\"{u}rlimann 2006) could be
readily generalized. Finally, it has not escaped our notice that
several applications recently advanced in the context of Benford's
law, such as fraud detection or stock market analysis (Nigrini
2000), could eventually be generalized to the wider context of GBL
formalism. This generalization also extends to stochastic sieve
theory (Hawkins 1957), dynamical systems that follow Benford's law
(Berger  \emph{et al.} 2005, Miller \emph{et al.} 2006) and their
relation to stochastic
multiplicative processes (Manrubia \emph{et al.} 1999).\\

\begin{acknowledgements}
We thank I. Parra for helpful suggestions and O. Miramontes, J.
Bascompte, D.H. Zanette and S.C. Manrubia for comments on a previous
draft. This work was supported by grant FIS2006-08607 from the
Spanish Ministry of Science.
\end{acknowledgements}

\appendix{Statistical methods and technical digressions}
\subsection{\textbf{How to pick an integer \emph{at random}?}}
\subsubsection{Visualizing the Generalized Benford law pattern in
prime numbers as a biased random walk}
\begin{figure}[h]
\centering
\includegraphics[width=0.7\textwidth]{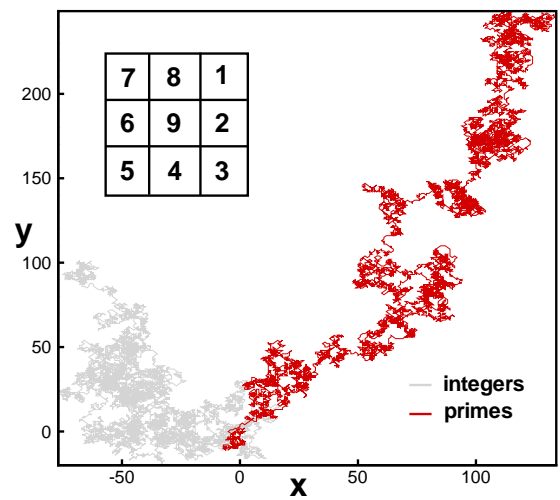}
\hspace{2cm} \caption{\small{Random walks. Grey: 2D Random walk in
which at each step we pick at random a natural from $[1,10^6]$ and
move forward depending on the value of its first significative digit
following the rules depicted in the inner table. The behavior
approximates an uncorrelated Brownian motion: integers first digit
is uniformly distributed. Red: same random walk but picking at
random primes in $[1,10^6]$: in this case the random walk is clearly
biased.}}\label{prw}
\end{figure}

 \noindent In order the
pattern already captured in figure 1 of the main text to become more
evident, we have built the following 2D random walk
\begin{eqnarray}
x(t+1)=x(t)+\xi_x\nonumber\\
y(t+1)=y(t)+\xi_y,
\end{eqnarray}
where $x$ and $y$ are cartesian variables with $x(0)=y(0)=0$, and
both $\xi_x$ and $\xi_y$ are discrete random variables that take
values $\in \{0,-1,1\}$ depending on the first digit $d$ of the
numbers randomly chosen at each time step, according to the rules
depicted in figure \ref{prw}. Thereby, in each iteration we peak at
random a positive integer (grey random walk) or a prime (red random
walk) from the interval $[1,10^6]$, and depending on its first
significative digit $d$, the random walker moves accordingly (for
instance if we peak prime $13$, we have $d=1$ and the random walker
rules provide $\xi_x=1$ and $\xi_y=1$: the random walker moves
up-right). We have plotted the results of this 2D Random Walk in
figure \ref{prw} for random picking of integers (grey random walk)
and for random picking of primes (red random walk). Note that while
the grey random walk seems to be a typical uncorrelated Brownian
motion (enhancing the fact that the first digit distribution of the
integers is uniformly distributed), the red random walk is clearly
biased: this is indeed a visual characterization of the pattern.
Observe that if the interval in which we randomly peak either the
integers or the primes wasn't of the shape [$1,10^D$], there would
be a systematic bias present in the pool and consequently both
integer and prime random walks would be biased: it comes thus
necessary to define the intervals under study in that way.

\subsubsection{Natural density} \noindent If primes were for instance
Benford distributed, one should expect that if we pick a prime at
random, this one should start by number 1 around $30\%$ of the time.
But what does the sentence '\emph{Pick a prime at random}' stand
for? Notice that in the previous experiment (the 2D biased Random
Walk) we have drawn whether integers or primes at random from the
pool $[1,10^6]$. All over the paper, the intervals $[1,N]$ have been
chosen such that $N=10^D$, $D \in \mathbb{N}$. This choice isn't
arbitrary, much on the contrary, it relies on the fact that whenever
studying infinite integer sequences, the results strongly depend on
the interval under study. For instance, everyone will agree that
intuitively the set of positive integers $\mathbb{N}$ is an infinite
sequence whose first digit is uniformly distributed: there exist as
many naturals starting by one as naturals starting by nine. However
there exist subtle difficulties at this point that come from the
fact that the first digit natural density is not well defined. Since
there exist infinite integers in $\mathbb{N}$ and consequently it is
not straightforward to quantify the quote 'pick an integer at
random' in a way in which satisfies the laws of probability, in
order to check if integers have a uniform distributed first
significant digit, we have to consider finite intervals $[1,N]$.
Hereafter, notice that uniformity \emph{a priori} is only respected
when $N=10^D$. For instance, if we choose the interval to be
$[1,2000]$, a random drawing from this interval will be a number
starting by $1$ with high probability, as there are obviously more
numbers starting by one in that interval. If we increase the
interval to say $[1,3000]$, then the probability of drawing a number
starting by $1$ or $2$ will be larger than any other. We can easily
come to the conclusion that the first digit density will oscillate
repeatedly by decades as $N$ increases without reaching convergence,
and it is thereby said that the set of positive integers with
leading digit $d$ ($d=1,2...,9$) does not possess a natural density
among the integers. Note that the same phenomenon is likely to take
place for the primes (see Chris Caldwell's \emph{The Prime Pages}
for an introductory discussion
in natural density and Benford's law for prime numbers and references therein).\\

\noindent In order to overcome this subtle point, one can: (i)
choose intervals of the shape $[1,10^D]$, where every leading digit
has equal probability a priori of being picked. According to this
situation, positive integers $\mathbb{N}$ have a uniform first digit
distribution, and in this sense Diaconis (1977) showed that primes
do not obey Benford's law as their first digit distribution is
asymptotically uniform. Or (ii) use average and summability methods
such as the Cesaro or the logarithm matrix method $\ell$ (Raimi
1976) in order to define a proper first digit density that holds in
the infinite limit. Some authors have shown that in this case, both
the primes and the integers are said to be \emph{weak} Benford
sequences (Raimi 1976, Flehinger 1966, Whitney 1972).\\

\noindent As we are dealing with finite subsets and in order to
check if a pattern \emph{really} takes place for the primes, in this
work we have chosen intervals of the shape $[1,10^D]$ to assure that
samples are unbiased and
that all first digits are equiprobable \emph{a priori}. \\

\subsection{\textbf{Statistical methods}}
\subsubsection{Method of moments} \noindent In order to estimate the
best fitting between a GBL with parameter $\alpha$ and a data set,
we have employed the method-of-moments. If GBL fits the empirical
data, then both distributions have the same first moments, and the
following relation holds:
\begin{equation}
\sum_{d=1}^9 {dP(d)}=\sum_{d=1}^9 {d P^e(d)}\label{momentos}
\end{equation}
where $P(d)$ and $P^e(d)$ are the observed normalized frequencies
and GB expected probabilities for digit $d$, respectively. Using a
Newton-Raphson method and
iterating equation \ref{momentos} until convergence, we have calculated $\alpha$ for each sample $[1,N]$.\\

\subsubsection{Statistical tests} \noindent Typically, chi-square
goodness-of-fit test has been used in association with Benford's Law
(Nigrini 2000). Our null hypothesis here is that the sequence of
primes follow a GBL. The test statistic is:
\begin{equation}
\chi^2 = M \sum_{d=1}^9 {{\big(P(d)- P^e(d)\big)^2} \over{P^e(d)}},
\label{chi}
\end{equation}
where $M$ denotes the number of primes in $[1,N]$. Since we are
computing parameter $\alpha(N)$ using the mean of the distribution,
the test statistic follows a $\chi^2$ distribution with $9-2=7$
degrees of freedom, so the null hypothesis is rejected if
$\chi^2>\chi^2_{a,7}$, where $a$ is the level of significance. The
critical values for the $10\%$, $5\%$, and $1\%$ are $12.02$,
$14.07$, and $18.47$ respectively. As we can see in table
\ref{Table1}, despite the excellent visual agreement (figure 1 in
the main text), the $\chi^2$ statistic goes up with sample size and
consequently the null hypothesis can't be rejected only for
relatively small sample sizes $N<10^{9}$. As a matter of fact,
chi-square statistic suffers from the excess power problem on the
basis that it is size sensitive: for huge data sets, even quite
small differences are statistically significant (Nigrini 2000). A
second alternative is to use the standard $Z$-statistics to test
significant differences. However, this test is also size dependent,
and hence registers the same problems as $\chi^2$ for large samples.
Due to this facts, Nigrini (2000) recommends for Benford analysis a
distance measure test called Mean Absolute Deviation (MAD). This
test computes the average of the nine absolute differences between
the empirical proportions of a digit and the ones expected by the
GBL. That is:
\begin{equation}
\textrm{MAD} = {1 \over 9} \sum_{d=1}^9 {\big|P(d)- P^e(d)\big|}
\label{MAD}
\end{equation}
This test overcomes the excess power problem of $\chi^2$ as long as
it is not influenced by the size of the data set. While MAD lacks of
cut-off level, Nigrini (2000) suggests that the guidelines for
measuring conformity of the first digits to Benford Law to be: MAD
between $0$ and $0.4\cdot10^{-2}$ imply close conformity, from
$0.4\cdot10^{-2}$ to $0.8\cdot10^{-2}$ acceptable conformity, from
$0.8\cdot10^{-2}$ to $0.12\cdot10^{-1}$ marginally acceptable
conformity, and finally, greater than $0.12\cdot10^{-1}$,
nonconformity. Under these cut-off levels we can not reject the
hypothesis that the first digit frequency of the prime numbers
sequence follows a GBL. In addition,
the Maximum Absolute Deviation $m$ defined as the largest term of MAD is also showed in each case.\\

\noindent As a final approach to testing for a similarity between
the two histograms, we can check the correlation between the
empirical and theoretical proportions by the simple regression
correlation coefficient $r$ in a scatterplot. As we can see in table
\ref{Table1} the empirical data are highly correlated with a
Generalized Benford distribution.\\

\noindent The same statistical tests have been performed for the
case of the Riemann non trivial zeta zeros sequence (table
\ref{Table2}), with similar results.

\begin{table}
\center
\begin{tabular}{|c|c|c|c|c|c|c|}
\hline
$N$ & $M$ = \# primes & $\chi^2$ & $m$ & MAD &$r$\\
\hline
$10^4$ & $1229$ & $0.45$ &    $0.32\cdot10^{-2}$ & $0.19\cdot10^{-2}$ & $0.96965$\\
$10^5$ & $9592$ & $0.62$ &    $0.21\cdot10^{-2}$ & $0.65\cdot10^{-3}$ & $0.99053$\\
$10^6$ & $78498$ & $0.61$ &    $0.50\cdot10^{-3}$ & $0.26\cdot10^{-3}$ & $0.99826$\\
$10^7$ & $664579$ & $0.77$ &    $0.17\cdot10^{-3}$ & $0.11\cdot10^{-3}$ & $0.99964$\\
$10^8$ & $5761455$ & $2.2$ &    $0.15\cdot10^{-3}$ & $0.56\cdot10^{-4}$ & $0.99984$\\
$10^9$ & $50847534$ & $11.0$ &    $0.11\cdot10^{-3}$ & $0.42\cdot10^{-4}$ & $0.99988$\\
$10^{10}$ & $455052511$ & $61.2$ & $0.90\cdot10^{-4}$ & $0.33\cdot10^{-4}$ & $0.99991$\\
$10^{11}$ & $4118054813$ & $358.5$ & $0.74\cdot10^{-4}$ & $0.27\cdot10^{-4}$ & $0.99993$\\
\hline
\end{tabular}
\caption{\small{Table gathering the values of the following
statistics: $\chi^2$, Maximum Absolute Deviation ($m$), Mean
Absolute Deviation (MAD) and correlation coefficient ($r$) between
the observed first significant digit frequency of the set of $M$
primes in
 $[1,N]$ and the expected Generalized Benford distribution (eq.
 \ref{benfordgen})
 with an exponent $\alpha(N)$ given by eq. \ref{alfa_n} with $a=1.1$).
 While $\chi^2$-test rejects the hypothesis for very large
 samples due to its size sensitivity, every other test cannot reject it,
 enhancing the goodness-of-fit between the data and the GB
 distribution.}}
\label{Table1}
\end{table}

\begin{table}
\center
\begin{tabular}{|c|c|c|c|c|c|c|}
\hline
$N$ & $M$ = \# zeros & $\chi^2$ & $m$ & MAD &$r$\\
\hline

$10^3$ & $649$      &$0.14$& $0.32\cdot10^{-2}$ &    $0.13\cdot10^{-2}$ & $0.99701$\\
$10^4$ & $10142$    &$0.23$& $0.11\cdot10^{-2}$ &    $0.41\cdot10^{-3}$ & $0.99943$\\
$10^5$ & $138069$   &$0.75$& $0.54\cdot10^{-3}$ &    $0.20\cdot10^{-3}$ & $0.99974$\\
$10^6$ & $1747146$  &$3.6$& $0.34\cdot10^{-3}$ &    $0.13\cdot10^{-3}$ & $0.99983$\\
$10^7$ & $21136126$ &$20.3$& $0.23\cdot10^{-3}$ &    $0.86\cdot10^{-4}$ & $0.99988$\\
\hline
\end{tabular}
\caption{\small{Table gathering the values of the following
statistics: $\chi^2$, Maximum Absolute Deviation ($m$), Mean
Absolute Deviation (MAD) and correlation coefficient ($r$) between
the observed first significant digit frequency in the $M$ zeros in
$[0,N]$ and the expected Generalized Benford distribution (eq.
\ref{gbl2} with and exponent $\alpha(N)$ given by eq. \ref{alfa_n}
with $a=2.92$). While $\chi^2$-test rejects the hypothesis for very
large samples due to its size sensitivity, every other test can't
reject it,
 enhancing the goodness-of-fit between the data and the GB
 distribution.}}
 \label{Table2}
\end{table}

\subsection{\textbf{Cram\'{e}r's model}}
\noindent The prime number distribution is deterministic in the
sense that primes are determined by precise arithmetic rules.
However, its apparent local randomness has suggested several
stochastic interpretations. Concretely, Cram\'{e}r (1935, see also
Tenembaum 2000) defined the following model: assume that we have a
sequence of urns $U(n)$ where $n = 1,2,...$ and put black and white
balls in each urn such that the probability of drawing a white ball
in the $k^{th}$-urn goes like $1/\log k$. Then, in order to generate
a
 sequence of pseudo-random prime numbers we only need to draw a ball
from each urn: if the drawing from the $k^{th}$-urn is white, then
$k$ will be labeled as a pseudo-random prime. The prime number
sequence can indeed be understood as a concrete realization of this
stochastic process. With such model, Cram\'{e}r studied amongst
others the distribution of gaps between primes and the distribution
of twin primes as far as statistically speaking, these distributions
should be similar to the pseudo-random ones generated by his model.
Quoting Cram\'{e}r: `With respect to the ordinary prime numbers, it
is well known that, roughly speaking, we may say that the chance
that a given integer $n$ should be a prime is approximately $1/\log
n$. This suggests that by considering the following series of
independent trials we should obtain sequences of integers presenting
a certain analogy
with the sequence of ordinary prime numbers $p_n$'.\\

\noindent In this work we have simulated a Cram\'{e}r process, in
order to obtain a sample of pseudo-random primes in $[1,10^{11}]$.
Then, the same statistics performed for the prime number sequence
have been realized in this sample. Results are summarized in table
\ref{Tablecram1}. We can observe that the Cram\'{e}r's model
reproduces the same behavior, namely: (i) The first digit
distribution of the pseudo-random prime sequence follows a GBL with
a size-dependent exponent that follows eq. \ref{alfa_n}. (ii) The
number of pseudo-primes found in each decade matches statistically
speaking to the actual number of primes. (iii) The $\chi^2$-test
evidences the same problems of power for large data sets. Having in
mind that the sample elements in this model are independent (what is
not the case in the actual prime sequence), we can confirm that the
rejection of the null hypothesis by the $\chi^2$-test for huge data
sets is not related to a lack of data independence but much likely
to the test's size sensitivity. (iv) The rest of statistical
analysis is similar to the one previously performed in the prime number sequence.\\

\begin{table}
\center
\begin{tabular}{|c|c|c|c|c|c|c|}
\hline
$N$ & $M$ = \# pseudo-random primes & $\chi^2$ & $m$ & MAD &$r$\\
\hline
$10^4$ & $1189$ & $1.20$ &    $0.17\cdot10^{-1}$ & $0.92\cdot10^{-2}$ & $0.639577$\\
$10^5$ & $9673$ & $0.43$ &    $0.33\cdot10^{-2}$ & $0.21\cdot10^{-2}$ & $0.969031$\\
$10^6$ & $78693$ & $0.39$ &    $0.59\cdot10^{-3}$ & $0.14\cdot10^{-2}$ & $0.990322$\\
$10^7$ & $664894$ & $0.09$ &    $0.23\cdot10^{-3}$ & $0.99\cdot10^{-4}$ & $0.999626$\\
$10^8$ & $5762288$ & $0.24$ &    $0.15\cdot10^{-3}$ & $0.53\cdot10^{-4}$ & $0.999855$\\
$10^9$ & $50850064$ & $1.23$ &    $0.11\cdot10^{-3}$ & $0.42\cdot10^{-4}$ & $0.999892$\\
$10^{10}$ & $455062569$ & $6.84$ & $0.90\cdot10^{-4}$ & $0.33\cdot10^{-4}$ & $0.999914$\\
$10^{11}$ & $4118136330$ & $41.0$ & $0.73\cdot10^{-4}$ & $0.27\cdot10^{-4}$ & $0.999937$\\
\hline
\end{tabular}
\caption{\small{Table gathering the values of the following
statistics: $\chi^2$, Maximum Absolute Deviation ($m$), Mean
Absolute Deviation (MAD) and correlation coefficient ($r$) between
the observed first significant digit frequency in the Cram\'{e}r
model for $M$ pseudo-random primes in $[1,N]$ and the expected
Generalized Benford distribution (eq. \ref{benfordgen}
 with an exponent $\alpha(N)$ given by eq. \ref{alfa_n} with $a=1.1$).}} \label{Tablecram1}
\end{table}

\smallskip

\newpage
\newpage
\newpage
\newpage
\newpage

\end{document}